\documentclass[12pt]{amsart}
\usepackage{amscd,amssymb}
\usepackage[arrow,matrix]{xy}
\usepackage[plainpages,backref,urlcolor=blue]{hyperref}

\theoremstyle{plain}
\newtheorem{thm}[subsection]{Theorem}

\newtheorem{prop}[subsection]{Proposition}
\newtheorem{cor}[subsection]{Corollary}

\theoremstyle{definition}
\newtheorem{rk}[subsection]{Remark}
\newtheorem{definition}[subsection]{Definition}
\newtheorem{ex}[subsection]{Example}

\numberwithin{equation}{section}
\newcommand{\Z}{\mathbb{Z}}
\newcommand{\Q}{\mathbb{Q}}

\newcommand{\C}{\mathbb{C}}

\newcommand{\PP}{\mathbb{P}}
\newcommand{\F}{\mathbb{F}}

\DeclareMathOperator{\Tors}{Tors}

\begin{document}

\title [On the topology of some quasi-projective surfaces]
{On the topology of some quasi-projective surfaces}

\author[Alexandru Dimca]{Alexandru Dimca$^1$}
\address{Univ. Nice Sophia Antipolis, CNRS,  LJAD, UMR 7351, 06100 Nice, France. }
\email{dimca@unice.fr}

\thanks{$^1$ Partially supported by Institut Universitaire de France}



\begin{abstract} Let $X$ be surface with isolated singularities in the complex projective space $\PP^3$ and let denote $Y$ the smooth part of $X$. In this note we discuss, mostly on specific examples, some aspects of the topology of such quasi-projective surfaces $Y$: the fundamental groups and the associated Galois coverings, the second homotopy groups and the mixed Hodge structure on the first cohomology group.

\end{abstract}

\maketitle


\section{Introduction and statements of results } \label{sec:intro}

Surfaces theory is a classical subject, with a very rich history and a number of excellent textbooks,
as for instance \cite{Ba}, \cite{Be}. It is quite surpring that new open questions arise even in such a classical subject, and one of the purposes of this note is to state such an open question in Example \ref{ex1} related to the cubic surfaces in $\PP^3$. These surfaces have been classified already by Cayley, see for a modern presentation \cite{BW}. Another open question in relation to the Zariski sextic with 6 cusps appears in Example \ref{ex2}.

Let $X$ be surface with isolated singularities in the complex projective space $\PP^3$. Then it is well known that $X$ is simply-connected, see for instance \cite{D2}. Let $Y$ denote the smooth part of $X$. So $Y$ is obtained from $X$ be removing a finite number of points. However, unlike the case when $X$ is smooth, this operation alters sometimes the fundamental groups and we may get 
quasi-projective surfaces $Y$ with $\pi_1(Y) \ne 0$. We omit the base points in this note, since our spaces $Y$ are path-connected, hence the isomorphism class of $\pi_1(Y,y)$ is independent of $y \in Y$.
For related results on fundamental groups of surfaces we refer to \cite{ADH}, \cite{DPS1} and \cite{DPS2}.

The first result describes the first integral homology group $H_1(Y)$ of the surface $Y$, which is exactly  the abelianization of the fundamental group $\pi_1(Y)$.
In the sequel $\Tors$ denotes the torsion part of a finitely generated abelian group.
\begin{thm}
\label{prop1}
Let $X$ be a surface with isolated singularities in $\PP^3$, let $Z$ be the singular set of $X$ and set $Y=X \setminus Z$. Then one has the following.

\noindent (i) $H_1(Y)=  H^3(X).$
In particular, if $X$ is a $\Q$-manifold, e.g. when $X$ has only simple singularities of type $A_n$, $D_n$, $E_6$, $E_7$ and $E_8$, then $H_1(Y)= \Tors H_2(X)$ is a finite group.

\noindent (ii) $H_3(Y)$ is a free abelian group of rank given by $|Z|-1$. 
\end{thm}
Note that the surface $Y$ is never affine, since a regular function $\phi$ defined on $Y$ extends to $X$, as $X$ is normal, and hence $\phi$ has to be constant.

The next result shows that there is a geometrically induced epimorphism 
$$\Gamma_g \to \pi_1(Y),$$
 where $\Gamma_g$ is the fundamental group of a smooth plane curve of genus
$$g =\frac{(d-1)(d-2)}{2}.$$
\begin{prop}
\label{prop2}
(i) For a generic plane $H$ in $\PP^2$, the intersection $C=X \cap H$ is a smooth curve contained in $Y$ and the inclusion $i:C \to Y$ induces an epimorphism
$$i_{\sharp}: \pi_1(C) \to  \pi_1(Y).$$
(ii) For any plane $H$ in $\PP^2$ such that the intersection $C=X \cap H$ is a (possibly singular) curve contained in $Y$, the inclusion $i:C \to Y$ induces an epimorphism
$$i_{\sharp}: \pi_1(C) \to  \pi_1(Y).$$
In particular, if $Y$ contains a rational cuspidal plane curve $C$, then $\pi_1(Y)=0$.
\end{prop}

\begin{cor}
\label{c1} Let $X$ be a surface with isolated singularities in $\PP^3$.
If $X$ is a surface of degree 3, then the fundamental group  $\pi_1(Y)$ of its smooth part $Y$ is abelian. More precisely, denote by $X(4A_1)$ the cubic surface in $\PP^3$ having as singularities 4 nodes $A_1$, and similarly for $X(3A_2)$, $X(A_1A_5)$ and $X(2A_1A_3)$. Then the corresponding smooth quasi-projective surfaces $Y(4A_1)$, $Y(3A_2)$, $Y(A_1A_5)$ and $Y(2A_1A_3)$ have the following fundamental groups.
$$ \pi_1(Y(4A_1))=\pi_1(Y(A_1A_5))=\pi_1(Y(2A_1A_3))=\Z/2\Z \text{ and } \pi_1(Y(3A_2))=\Z/3\Z.$$
\end{cor}

For a description of the associated Galois coverings and the second homotopy groups of these surfaces see Example \ref{ex1}.
\begin{rk}
\label{rk0.1}
If $X$ is a cubic surface with simple singularities in $\PP^3$, then $X$ is an exemple of a log del Pezzo surface, see \cite{AN} for the corresponding definition. It is known that the fundamental group of the smooth part $Y$ of such a log del Pezzo surface is always finite, see  \cite{GZ}, \cite{KMc}, \cite{Xu1}. For an extension of this result to higher dimensional log Fano varieties see
\cite{Keb}, \cite{Xu2}.

\end{rk}

We have also the following.

\begin{prop}
\label{prop2.5}
Let $X$ be a surface with isolated singularities in $\PP^3$, let $Z$ be the singular set of $X$ and set $Y=X \setminus Z$. For each singular point $z \in Z$, let $L_z$ denote the link of the singularity $(X,z)$. Then there is morphism
$$\Pi_{z\in Z} \pi_1(L_z) \to \pi_1(Y),$$
where $\Pi_{z\in Z} \pi_1(L_z)$ denotes the free product of the family of groups  $(\pi_1(L_z))_{z\in Z}$, whose image is not contained in any normal subgroup of $\pi_1( Y)$.
\end{prop}

One can ask about the mixed Hodge structures on the surfaces $Y$. Here is the answer.

\begin{thm}
\label{thm2}
Let $X$ be surface with isolated singularities in $\PP^3$, let $Z$ be the singular set of $X$ and set $Y=X \setminus Z$. Then  $H_c^3(Y)=  H^3(X)$ is a pure Hodge structure of weight 3 and by duality, 
 $H^1(Y)$ is a pure Hodge structure of weight 1. 

\end{thm}

The cohomology groups $H_c^2(Y)=  H^2(X)$ have a more complicated mixed Hodge structure, in general with several possible weights, for more on this subject see \cite{DHam} and \cite{DL}.

The Hodge theoretic result  in Theorem \ref{thm2} has the following pure topological consequence. Let $C:g(x,y,z)=0$ be a reduced  curve in $\PP^2$ of degree $d$. Consider the surfaces $X_C: g(x,y,z)+t^d=0$, whose singularities are the degree $d$ suspension of the singularities of the curve $C$. Let $Y_C$ be the smooth part of $X_C$ as above and denote by $F_C$ the Milnor fiber of $g$, namely the affine smooth surface $$F_C: \  \  g(x,y,z)=1$$
in $\C^3$. Then clearly $F_C \subset Y_C$ is a Zariski open subset and hence the inclusion $j:F_C \to Y_C$ induces an epimorphism $\pi_1(F_C) \to \pi_1(Y_C)$. By duality, we get a monomorphism $j^*:H^1(Y_C, \Q) \to H^1(F_C,\Q)$.

\begin{prop}
\label{prop3}
The image of the monomorphism $j^*:H^1(Y_C, \Q) \to H^1(F_C,\Q)$ is exactly 
$H^1(F_C,\Q)_{\ne 1}$, the non-fixed part of $H^1(F_C,\Q)$ under the monodromy action.
\end{prop}
This result allows us to construct many examples of surfaces $X_C$ such that $H^1(Y_C, \Q)$ (and presumably $\pi_1(Y)$) is quite large, e.g. using as $C$ various line arrangements in $\PP^2$, see \cite{Su},   \cite{Su2}, \cite{PStrip} for various monodromy computations in this case. To increase these groups, one may also use non-linear arrangements as well, as described for instance in Example 5.14 in \cite{D3}. One example using Zariski sextic curve with 6 cusps is given in Example \ref{ex2} below.

\medskip

Some open questions appear in Example \ref{ex1} and Example \ref{ex2}. A major open question is to develop a general strategy for the computation of the fundamental groups for this class of surfaces.

\medskip

This note gives a number of (very limited) answers to questions that Ciro Ciliberto asked me some time ago. I would like to thank him for asking the questions and for explaining to me the role of del Pezzo surfaces in Example \ref{ex1} and Remark \ref{rk3}. Many thanks also to 
De-Qi Zhang for his comments on the fundamental groups of log del Pezzo surfaces and open $K3$ surfaces, incorporated essentially in Remarks \ref{rk0.1} and  \ref{rk4} below. 

\section{The proofs and additional examples} 

We consider first Theorem \ref{prop1}. We apply Lefschetz duality theorem, see for instance \cite{Sp}, p. 297 to the compact relative 4-manifold $(X,Z)$, where $Z$ is the finite set of singular points of the surface $X$. Since any pair of algebraic sets is triangulable, it follows that the pair 
$(X,Z)$ is taut, see \cite{Sp}, p. 291, and hence we have an isomorphism
$H_k(Y)=H^{4-k}(X,Z)$.  To prove the first claim (i), we take $k=1$ and 
the long exact sequence of the pair $(X,Z)$ yields an isomorphism $H^3(X,Z)=H^3(X)$. If $X$ is a $\Q$-manifold, it follows that $b_3(X)=b_1(X)=0$ and hence $H^3(X)$ is a finite group. It remains to use the standard fact that $\Tors H_2(X)=\Tors H^3(X)$, see \cite{Sp}, p. 244. For the classification of simple singularities $A,D,E$ and their properties we refer to \cite{D1}, \cite{D2}.
To prove the claim (ii), we take $k=3$ in the above isomorphism and get $H_3(Y)=H^1(X,Z)$. Then the 
long exact sequence of the pair $(X,Z)$ yields 
$$b_3(Y)= |Z|-1.$$

The proof of Proposition \ref{prop2} follows from the Zariski theorem of Lefschetz type stated for instance in \cite{D2}, p. 25 and the fact that $X$ admits the obvious Whitney regular stratification given by $Y$ and the finite set $Z$, see for instance in \cite{D2}, p. 5. For the part (ii), one has to use the careful description of "good" hyperplanes in this case given in loc.cit.
Moreover, an irreducible projective curve is simply-connected if and only if it is a rational cuspidal curve, i.e. $C$ is rational and any singular point of $C$ is unibranch.

To prove Corollary \ref{c1}, for the first claim we just apply Proposition \ref{prop2} and the fact that $\Gamma_1=\Z^2$ to get that $\pi_1(Y)$ is abelian. 

For the computation of the fundamental groups of the surfaces,
$Y(4A_1)$, $Y(3A_2)$, $Y(A_1A_5)$ and $Y(2A_1A_3)$, we use the corresponding results for the second homology groups described in \cite{D2}, p. 165.

The proof of Proposition \ref{prop2.5} follows by applying the van Kampen theorem to the open covering $Y,Z'$ of $X$, where $Z'$ is obtained as follows. Take a mimal tree $T$ formed of simple non-intersecting arcs connecting the points in $Z$. Add for each $z \in Z$ a small contractible neighborhood $B_z$ of $z$ in $X$ such that $B_z^*=B_z \setminus \{z\}$ is homotopically equivalent to the link $L_z$. Then $Z'$ is a tubular open neighborhood of $T \cup \cup_{z\in Z}B_z$. It follows that $Y \cap Z'$ has the homotopy type of the join of the links $L_z$, and hence 
$$\pi_1(Y \cap Z')=\Pi_{z\in Z} \pi_1(L_z) .$$
On the other hand, $Z'$ is contractible and $X$ is simply-connected, so van Kampen theorem implies that the inclusion $ Y \cap Z' \to Y$ induces a morphism $ \pi_1(Y \cap Z' )\to \pi_1( Y)$
whose image is not contained in any normal subgroup of $\pi_1( Y)$.

We consider now Theorem \ref{thm2}. Note that the first part of this proof gives an alternative proof for the claim (i) in Theorem \ref{prop1}. The exact sequence with compact supports for the pair $(X,Z)$ yields the isomorphism $H^3_c(Y)=H^3(X)$ of MHS (short for mixed Hodge structures).
The result follows using the fact that $H^3(X)$ is a pure Hodge structure, see \cite{St}. For duality between $H^3_c(Y)$ and $H^1(Y)$ we refer to \cite{PS}. 

Finally, to prove Proposition \ref{prop3}, we recall that we have a splitting
$$H^1(F_C,\Q)=H^1(F_C,\Q)_{1}\oplus H^1(F_C,\Q)_{\ne 1},$$
where $H^1(F_C,\Q)_{1}$ has pure weight 2 and $H^1(F_C,\Q)_{\ne 1}$ has pure weight 1, see \cite{DL}. Moreover, it is shown in \cite{DL} and in \cite{DP} that 
$$\dim H^1(Y_C, \Q)=\dim H^3(X_C, \Q)=\dim H^1(F_C,\Q)_{\ne 1}.$$
This clearly completes the proof.

\begin{rk}
\label{rk1}

It follows from Theorem \ref{prop1} (ii) that the three surfaces $X(4A_1)$, $X(A_1A_5)$ and $X(2A_1A_3)$ have distinct third Betti numbers, namely $3$, $1$ and $2$, hence they are not homotopically equivalent to each other. 
\end{rk}

\begin{ex}
\label{ex1} 
It follows from Corollary \ref{c1}, that each of  the surfaces,
$Y(4A_1)$, $Y(3A_2)$, $Y(A_1A_5)$ and $Y(2A_1A_3)$ has a finite unramified cover which is a simply-connected surface, i.e. the corresponding universal covering. In the case of the surface $Y(3A_2)$, we can chose the equation for
$X(3A_2)$ to be 
$$f=xyz-t^3=0$$
and hence the map $p:\PP^2 \to X(3A_2)$ given by
$$p([u,v,w])=[u^3,v^3,w^3,uvw]$$
is ramified precisely over the singular set $Z$, see also \cite{D2}, p. 166. Hence the universal cover of $Y(3A_2)$ is obtained from $\PP^2$ by deleting 3 points. This also implies $\pi_2(Y(3A_2))= \pi_2(\PP^2)= \Z$.

For the other three universal covering surfaces, the construction involves Cremona transformations and del Pezzo surfaces of degree 6, i.e. surfaces obtained from $\PP^2$ by blowing up 3 points. The easiest case to describe is for the surface $X(4A_1)$. Consider the classical Cremena rational map
$$ c_1: \PP^2 \to \PP^2,  \  \  \  [x,y,z] \mapsto [yz,xz,xy].$$
This map has $p_1=[1,0,0]$, $p_2=[0,1,0]$ and $p_3=[0,0,1]$ as indeterminacy points, see Example 1.5.1 in \cite{Do1}. Let $S$ denote the del Pezzo surface obtained by blowing-up these 3 points and note that $c_1$ lifts to a regular map $c_1':S \to S$, which is an involution, i.e. $c_1'^2=Id$, and has 4 fixed points, namely the points in $S$ corresponding to the points $[1,\pm 1,\pm1]$ in $\PP^2$. It follows that the quotient surface $S/<1,c_1'>$, which has 4 nodes, can be identified to $X(4A_1)$. This implies $\pi_2(Y(4A_1))= \pi_2(S)= \Z^4$, by an easy application of Hurewicz Theorem.

To get the universal cover of $Y(2A_1A_3)$, we start with the {\it first  degenerate standard quadratic transformation} 
$$ c_2: \PP^2 \to \PP^2,  \  \  \  [x,y,z] \mapsto [y^2,xy,xz],$$
see \cite{Do1}, p. 15, which is also an involution, i.e. $c_2^2=Id$, and has $p_1$ and $p_3$ as indeterminacy points.
This map has 2 lines of fixed points, namely the lines $x \pm y=0$ meeting at the point $p_3$.
After blowing up $p_1$, $p_3$ and an infinitely near point $p_2'$, we get as above  the degree 2 covering over $X(2A_1A_3)$. The $A_3$ singularity occurs as the $\Z/(2)$-quotient of an $A_1$ singularity, which in turn comes from the contraction of a $(-2)$-curve obtained in the blowing-up process, see \cite{Do1}, p. 15.

Finally, to get the universal cover of $Y(A_1A_5)$, we start with the 
{\it second degenerate standard quadratic transformation }
$$ c_3: \PP^2 \to \PP^2,  \  \  \  [x,y,z] \mapsto [x^2,xy,y^2-xz],$$
see \cite{Do1}, p. 16, which is also an involution, i.e. $c_3^2=Id$, and has  $p_3$ as its unique indeterminacy point.
This map has a smooth conic of fixed points, namely  $y^2-2xz=0$, passing through the point $p_3$.
After blowing up $p_3$ and 2 infinitely near points $p_1'$ and $p_2'$, we get as above  the degree 2 covering over $X(A_1A_5)$. The $A_5$ singularity occurs as the $\Z/(2)$-quotient of an $A_2$ singularity, which in turn comes from the contraction of a chain of two $(-2)$-curves obtained in the blowing-up process, see \cite{Do1}, p. 16.
\end{ex}

\begin{ex}
\label{ex2} 
Here is an example of a surface $X$ such that the fundamental group $\pi_1(Y)$ is infinite.
Let $X$ be the surface given by
$$(x^2+y^2)^3+(y^3+z^3)^2+t^6=0,$$
the degree 6 cover of $\PP^2$ ramified over the Zariski sextic with 6 cusps on a conic.
It is well known that $b_3(X)=2$, see for instance in \cite{D2}, p. 210.
It would be interesting to check whether $H^3(X)=\Z^2$ and also to determine the 
fundamental group $\pi_1(Y)$ in this classical example. 
\end{ex}

\begin{rk}
\label{rk2}

When $b_1(Y)=b_3(X)>0$, then one can use Theorem 2 in \cite{DPS2} to compute the germs at the origin of the characteristic varieties $V^1_r(Y)_1$ in terms of the resonance varieties germs $R^1_r(\tilde X)_0$, where $\tilde X$ is a resolution of singularities for $X$. Note that there is no dominant morphism  from $Y$ to non proper curve $C$, since such a curve is necessarily affine and we can repeat the map extension described after Theorem \ref{prop1}.
\end{rk}

\begin{rk}
\label{rk3}
Some of the considerations in this paper can be applied to a complete intersection surface $X$ with isolated singularities in $\PP^n$ with $n \geq 4$. The case when $X$ is the intersection of two quadrics in $\PP^4$ is considered in Proposition 4.3 in \cite{DTors}, where the corresponding group $H_1(Y)$ is computed and found to be $\Z/2\Z$ in two cases, namely for $Y(4A_1)$ and for $Y(2A_1A_3)$, with an obvious notation. A generic hyperplane section of $X$ in this case is again a curve of genus one, hence we get as in Corollary \ref{c1} that one has
$$\pi_1(Y)=H_1(Y) = \Z/2\Z,$$
in these cases as well. It would be interesting to find an analog of Proposition \ref{prop3} in the case of complete intersections of codimension $>1$. The corresponding degree 2 universal covering spaces associated with  $Y(4A_1)$ and  $Y(2A_1A_3)$  have a simple geometrical description in terms of del Pezzo surfaces of degree 8. There are two types of such surfaces.

The first type is the product $S=\PP^1 \times \PP^1$, which has the involution $\iota: S \to S$
given by $\iota ([x,y],[u,v])=([y,x],[v,u])$, which has 4 fixed points, namely $([1,\pm 1], [1,\pm 1])$. It can be shown that the quotient $S/<1,\iota>$ is nothing else but the surface $X(4A1)$.
This also implies $\pi_2(Y(4A_1))=\pi_2(S)=\Z^2.$

The universal covering space of $Y(2A_1A_3)$ can be described as follows.
On the minimal resolution $X'$ of $X=X(2A_1A_3)$ we have the three
$(-2)$-curves $E_1+E_2+E_3$ ($E_2$ being the ''central'' one) in the
resolution of the $A_3$ singularity and two $(-2)$-curves $E_4,E_5$ resolving the two $A_1$
points. The double cover of $Y$ gives rise
to a double  cover of $X'$ branched along the sum of some of the curves
$E_i$, which must be divisible by 2 in $Pic(X')$. The
branch divisor is $E_1+E_3+E_4+E_5$, thus the double cover $X''$ possesses
two $(-1)$-curves $E'_4,E'_5$ over $E_4,E_5$, and a cycle $E'_1+E'_2+E'_3$ over
$E_1+E_2+E_3$, with $E'_1,E'_3$ being $(-1)$-curves and $E_2^{' 2}=-4$. By
contracting all the $(-1)$ curves we see that $E'_2$ becomes a $(-2)$-curve and in fact we get
in this way a Hirzebruch surface $\F_2$. If we contract this last $(-2)$-curve, we get the cone $S$ in $\PP^3$ over a smooth conic $Q$, say $Q: x^2-yz=0$. Then $S$ is a singular (log) del Pezzo surface of degree 8 and it is the double cover of $X$ ramified exactly over the singular points.
The corresponding involution of $S$ can be taken to be $[x,y,z,t] \mapsto [-x,y,z,-t]$, whose
 fixed points are precisely

\noindent (i)  $[0,0,0,1]$,  the vertex of the cone,
which  gives the $A_3$- singularity in the quotient $X$, and

\noindent (ii) the two points $[0,1,0,0]$ and $[0,0,1,0]$, which give the two $A_1$ singularities of $X$.
This proof also implies $\pi_2(Y(2A_1A_3))=\pi_2(Q)=\Z.$
\end{rk}

\begin{rk}
\label{rk4}
Corollary \ref{c1} and Remark \ref{rk3} describe the fundamental groups of the smooth part $Y$ of some singular del Pezzo surfaces $X$, i.e. surfaces $X$ with an ample anticanonical divisor $-K_X$. For the reader convenience, and for their beauty, we recall below some similar results on the fundamental groups in the case of singular $K3$-surfaces, i.e. surfaces with $K_X=0$.

\noindent (i) Let first $X$ be a $K3$-surface with 16 $A_1$-singularities, for instance the surface in $\PP^3$ given by the equation
$$x^4+y^4+z^4-(x^2+y^2+z^2)t^2-x^2y^2-x^2z^2-y^2z^2+1=0.$$
Such a surface $X$ can be obtained as the quotient of an Abelian surface $A$ under the involution $a \mapsto -a$, see \cite{N}. Since $A$ is a quotient $\C^2/\Lambda$, with $\Lambda$ a lattice of rank 4, we get a map $p: \C^2 \to X$, presenting $X$ as the quotient of $\C^2$ under the non-commutative subgroup $\tilde \Lambda$ spanned by $\Lambda$ (regarded as translations) and the involution $v \mapsto -v$ inside the group of affine transformations of the plane $\C^2$. It follows that $\C^2 \setminus p^{-1}(Z) \to Y$ is a universal covering for $Y$, where $Z$ is the singular part of $X$, with deck transformation group $\pi_1(Y)=\tilde \Lambda$. We also get $\pi_2(Y)= \pi_2(\C^2 \setminus B )=0$, since $B= p^{-1}(Z)$ is discrete, hence of real codimension 4.

\noindent (ii) Let now $X$ be a $K3$-surface with 8 cusps $A_2$ and no other singularities. Some of these surface can embedded in $\PP^3$ and the corresponding equations can be found in 
\cite{Bar2}. For such surfaces, one has either $\pi_1(Y)= (\Z/3\Z),$ or $\pi_1(Y)= (\Z/3\Z)^2,$ see \cite{KZ}, Table 1, case $p=3$, $c=8$.

\noindent (iii) Let now $X$ be a $K3$-surface with 9 cusps $A_2$ and no other singularities. None of these surface can embedded in $\PP^3$, see 
\cite{Bar2}. For such surfaces, one has an extension
$$0 \to \Lambda \to \pi_1(Y) \to \Z/3\Z \to 0$$
where  $\Lambda$ a lattice of rank 4 as above, see \cite{KZ}, Table 1, case $p=3$, $c=8$.
Moreover, one has as above $\pi_1(Y)=\tilde \Lambda$, where $\tilde \Lambda$ is the subgroup of  the group of affine transformations of the plane $\C^2$ generated by $ \Lambda$ and an order 3 linear automorphism of $\C^2$ preserving $ \Lambda$, see \cite{OZ}, Example 1 and \cite{SI}. An alternative approach is described in \cite{Bar1}, and one may wonder if the fundamental groups $\pi_1(Y)$ are the same for these two distinct examples. 

As above, we also get $\pi_2(Y)=0.$

For other low degree cuspidal surfaces $X$ in $\PP^3$ admiting degree 3 unramified covering over their smooth part $Y$, in particular a discussion of irreducible families of such surfaces, see \cite{BR}.
\end{rk}

\end{document}